\documentclass[12pt]{article} 
{  
\begin{document} 
\begin{center} 
{\mathversion{bold}
{\LARGE {\bf Some deformations of $U[sl(2)]$\\[2mm] and their representations}}}
\\[9mm] 
{\bf Nguyen Anh Ky}\\[2mm] 
Institute of Physics\\ 
P.O. Box 429, Bo Ho, Hanoi 10 000, Vietnam 
\end{center} 
\vspace*{2mm} 
\begin{abstract}     

Some one- and two-parametric deformations of  U[sl(2)] and their representations  
are considered. Interestingly, a newly introduced two-parametric deformation 
admits a class of infinite - dimensional representations which have no classical 
(non-deformed) and one-parametric deformation analogues, even at generic 
deformation parameters. 
\end{abstract} 
\begin{center}
\end{center}
\vspace*{5mm}
{\Large {\bf 1. Introduction}} 
 
  For over two decades, quantum  groups \cite{frt} - \cite{woro} have been a 
subject of great interest in physics and mathematics. 
A number of mathematical structures and physical applications of quantum  
groups have been obtained and investigated in details \cite{collect} - \cite{majid}.  
Depending on points of view, quantum groups  can be approached to in several ways. 
One of the approaches to quantum groups are the so-called Drinfel'd-Jimbo deformation 
of universal enveloping algebras \cite{drin, jim}. It is shown that the quantum groups 
of this kind (called also quantum  algebras) are non-commutative and 
non-cocommutative Hopf algebras \cite{drin}.  By construction, such a quantum group 
depends on one or  more parameters which are, in general, complex. One-parametric 
deformations have been well investigated and understood in various aspects,  while 
multi-parametric deformations are less understood, in spite of  some research progress 
made lately (see \cite{swz} - \cite{pq01} and references therein).  In this report, we 
consider some one-  and two-parametric deformations of the universal enveloping 
algebra $U[sl(2)]$ and their representations. 

  The quantum  group  $U_q[sl(2)]$ as an one-parametric deformation of the universal 
enveloping algebra $U[sl(2)]$ is one of the best  investigated quantum groups 
\cite{kure} - \cite{pas}. 
What about two-parametric deformations of $U[sl(2)]$, they have been considered 
in several versions and in different aspects (see, for example, \cite{swz} - \cite{pq96} 
and references therein) though some of them  are, in fact, equivalent to one-parametric 
deformations (upto some rescales).  In this report we consider only two of them, 
denoted here by $U^{(1)}_{pq}[sl(2)]$ and $U^{(2)}_{pq}[sl(2)]$ (a common 
notation is $U_{pq}[sl(2)]$).  The first quantum group $U^{(1)}_{pq}[sl(2)]$ was 
already given before as a subgroup of a quantum supergroup \cite{pq96,pq00,pq01}, 
while the second one, $U^{(2)}_{pq}[sl(2)]$, is being introduced now.  It turns out that 
this new quantum group $U^{(2)}_{pq}[sl(2)]$ admits  a class of infinite - dimensional 
representations which have no  analogue in either the case of the non-deformed $sl(2)$ 
or the cases of previously introduced one- and two-parametric deformations of $U[sl(2)]$. 
For consistency,  we first consider  the non-deformed $sl(2)$, the one-parametric 
deformation $U_q[sl(2)]$ and the two-parametric deformations $U^{(1)}_{pq}[sl(2)]$ 
in the next three sections, Sect. 2,  Sect. 3  and Sect. 4, respectively. Then the 
two-parametric deformations $U^{(2)}_{pq}[sl(2)]$ is introduced and considered in 
Sect. 5.  Some discussions and conclusion are made in the last section, Sect. 6. Let us 
start now with $sl(2)$.\\[5mm]
{\Large {\bf 2. {\mathversion{bold}$sl(2)$} and representations}} 

As is well known, $sl(2)$ can be generated by three generators, say $E_+$, $E_-$ 
and $H$,  satisfying the commutation relations 
$$[H,E_{\pm}]=\pm E_{\pm},  \qquad [E_+,E_-]=2H.\eqno(1)$$ 
Demanding $(H)^{\dagger}=H$ and $(E_{\pm})^{\dagger}=E_{\mp}$, a 
representtaion induced from  a (normalised) highest weight state $|j,j\rangle$, 
$$H|j,j\rangle=j|j,j\rangle, \quad E_+|j,j\rangle=0,\eqno(2)$$ 
has the matrix elements 
$$H~|j,m\rangle = m ~|j,m\rangle, $$ 
$$ E_\pm~| j,m\rangle = 
\sqrt{(j\mp m)(j\pm m +1)}~ |j, m\pm 1\rangle,\eqno(3)$$ 
where $|j,m\rangle$, $m \leq j$, is one of the orthonormalised states which 
can be obtained by acting an appropriate monomial of $E_-$ on the highest 
weight state $|j,j\rangle$, 
$$|j,m\rangle= A_n(E_-)^n~|j,j\rangle,\qquad m=j-n, \eqno(4)$$ 
with $A_n$ a normalised coefficient, which for a non-negative (half) integral 
highest weight $j$ has the form
$$A_n=\sqrt{(2j)!\  n!\over (2j-n)!}\  .\eqno(5)$$ \\ 
In this case, the representations constructed  are simultaneously highest weight 
and lowest weight, that is, 
 $$E_+|j,j\rangle=0, \quad  E_-|j,-j\rangle=0,\eqno(6)$$ 
so, they are finite-dimensional (and also unitary and irreducible) of dimension 
$2j+1$ (with  $m$ taking by steps of 1 all values in the range $-j\leq m \leq j$). 
The situation is similar in the case of the one-parametric quantum group $U_q[sl(2)]$ 
at a generic deformation parameter $q$ (i.e., at $q$ not being a root of unity). \\[5mm] 
{\Large {\bf 3.  One-parametric deformation {\mathversion{bold}$U_q[sl(2)]$}}}

   As mentioned above, $U_q[sl(2)]$ is one of the best investigated quantum groups. 
It can be defined through three generators $E_+$, $E_-$ and $H$  subject to the 
deformed defining relations 
$$[H,E_{\pm}]=\pm E_{\pm},  \quad [E_+,E_-]=[2H]_q,\eqno(7)$$
where 
$$[X]_q= {q^X-q^{-X}\over q-q^{-1}},\eqno(8)$$ 
is the so-called one-parametric deformation (or $q$-deformation, for short) of  an 
operator or a number $X$, with $q$ being a, complex in general, parameter.  At a 
generic $q$, the representation structure of $U_q[sl(2)]$ is  
similar to that of $sl(2)$, more precisely, a representation of $U_q[sl(2)]$ induced 
in the same way (as done for $sl(2)$) from a highest weight state $|j,j\rangle$ defined 
as in (2) has the matrix elements  $$H~|j,m\rangle = m ~|j,m\rangle, $$ 
$$ E^\pm~|j,m\rangle = 
\sqrt {[j\mp m]_q[j\pm m +1]_q}~ |j, m\pm 1\rangle.\eqno(9)$$\\ 
   These representations are again of finite dimension $2j+1$ (and irreducible as $q$ 
is generic) for non-negative (half) integral highest weights $j$, since Eqs. (6) are 
again satisfied (i.e., the representations are again highest weight and lowest weight 
simultaneously).\\[5mm] 
{\Large   {\bf 4. A two-parametric deformation {\mathversion{bold}$U^{(1)}_{pq}[sl(2)]$}}} 

The quantum group $U^{(1)}_{pq}[sl(2)]$ considered here is only one of the 
two-parametric deformations of $U[sl(2)]$ introduced in the literature.  
Appeared in \cite{pq96} as a subgroup of a quantum supergroup $U_{pq}[gl(2/1)]$, 
this quantum group  can be generated by three generators $E_+$, $E_-$ and $H$  
satisfying somewhat more twisted defining relations 
$$[H,E_{\pm}]=\pm E_{\pm},  \quad [E_+,E_-]= 
\left (p\over q\right )^{J-H}[2H]_{pq}, 
%\quad (H)^{\dagger}=H, \quad (E_{\pm})^{\dagger}=E_{\mp}, 
\eqno(10)$$ 
where $J$ is the so-called the "maximal spin" (or the hight weight) operator and 
$$[X]_{pq}= {q^X-p^{-X}\over q-p^{-1}},\eqno(11)$$ 
is a notation for a two-parametric deformation (or $pq$-deformation) of  an operator 
or a number $X$, with $p$ and $q$ being deformation parameters which are in 
general complex. Constructed in the same way for generic $p$ and $q$, the (highest weight) 
representations of $U^{(1)}_{pq}[sl(2)]$ are similar to those of $U_q[sl(2)]$ and $sl(2)$, 
$$H~|j,m\rangle = m ~|j,m\rangle,$$ 
$$ E^\pm~|j,m\rangle = 
\sqrt{[j\mp m]_{pq}[j\pm m +1]_{pq}}~ |j, m\pm 1\rangle.\eqno(12)$$ 
They are again finite - dimensional for non-negative (half) integral $j$.  Although the 
quantum supergroup $U_{pq}[gl(2/1)]$ is not equivalent to an one-parameter 
deformation of $U[gl(2/1)]$, its subgroup $U^{(1)}_{pq}[sl(2)]$ can be shown 
to be equivalent to an one-parametric deformation of $U[sl(2)]$. The situation  is, 
however, different in the next case.\\[5mm] 
{\Large {\bf 5. An alternative two-parametric deformation 
{\mathversion{bold}$U^{(2)}_{pq}[sl(2)]$}}} 

 Now we introduce an alternative two-parametric quantum  group $U^{(2)}_{pq}[sl(2)]$  
generated also by three generators  $E^+$, $E^-$ and $H$ but the latter now satisfy defining 
relations simpler than those of $U^{(1)}_{pq}[sl(2)]$,  
$$[H,E^{\pm}]=\pm E^{\pm},  \quad [E^+,E^-]=[H]_{pq} 
.\eqno(13)$$ 
Following the way of constructing highest weight representations in the previously 
considered cases,  we find  representations of $U^{(2)}_{pq}[sl(2)]$ corresponding to (3), 
(9) and (12),   
$$H~|j,m\rangle = m ~|j,m\rangle, $$ 
$$E^+~|j,m\rangle = 
\left ( ~ {q^{j+m+1} [j-m]_q-p^{-j-m-1} [j-m]_p}\over {q-p^{-1}}\right ) ^{1/2} 
~ |j, m+1\rangle,$$ 
$$\quad E^-~|j,m\rangle =  
\left  (~ {q^{j+m}[j-m+1]_q-p^{-j-m}[j-m+1]_p}\over {q-p^{-1}}\right ) ^{1/2} 
~ |j, m-1\rangle.\eqno(14)$$ \\ 
In comparison with the case of $U^{(1)}_{pq}[sl(2)]$, for the sake of simplicity of the  
difining relations, the matrix elements  in this case become a bit more complicated. 
Moreover, even at a non-negative (half) integral highest weights $j$, we observe 
that representations of the form (14) are, in general, infinite - dimensional, unlike 
those of  $sl(2)$,  $U_q[sl(2)]$ and $U^{(1)}_{pq}[sl(2)]$ in (3), (9) and (12), 
respectively. These 
representaions of $U^{(2)}_{pq}[sl(2)]$ are still highest weight by construction 
but they are no longer lowest weight for arbitrary $p$ and $q$. It is a new class of 
infinite - dimensional representations of $U_{pq}[sl(2)]$ not found before in the cases 
of $sl(2)$ and its previously considered deformations. Note that $U^{(2)}_{pq}[sl(2)]$ 
is by no means equivalent to the one-parametric $U_q[sl(2)]$
unless at $p=q$ and at, maybe, other special choices of $p$ and $q$.
\\[5mm] 
{\Large  {\bf 6. Conclusion}} 

   We have introduced a new quantum group $U^{(2)}_{pq}[sl(2)]$ which is a 
two-parametric deformation of $U[sl(2)]$. Interestingly, it is  showed that 
this quantum group admits a class of infinite - dimensional representations 
which have no analogues in the prviously considered cases of the non-deformed 
$sl(2)$, the one-parametric deformation $U_q[sl(2)]$ and other two-parametric 
deformations of $U[sl(2)]$.  It is a new phenomenon of $U^{(2)}_{pq}[sl(2)]$.\\[3mm]
{\bf Proposition 1}: {\it Highest weight representations of the two-parametric quantum 
algebra $U^{(2)}_{pq}[sl(2)]$ defined in (10) and (11) are in general infinite - dimensional, 
even for non-negative (half) integral highest weights}.
 
These representations are still unitary and, in general, irreducible. They may 
become reducible under certain additional conditions such as one or both of 
the deformation parameters to be roots of unity or  at some other special choices 
if allowed. 
However, in these cases, it is possible to extract finite - dimensional irreducible 
representations from the infinite - dimensional reducible ones. \\[4mm]
{\bf Proposition 2}: {\it The representations (14) are reducible and contain 
finite-dimensional subrepresentations iff the equation 
$$f(x)\equiv  q^{2j-x}[x+1]_q-p^{-2j+x}[x+1]_p=0\eqno(15)$$
has positive intergral solutions.}

Other classes of 
infinite - dimensional representations of $U_{pq}[sl(2)]$ may be found by using 
the method of 
\cite{stoyanov}. 
\\[4mm]  
{\bf Acknowledgement}: I would like to thank Randjbar - Daemi for kind hospitality 
at the Abdus Salam International Centre for Theoretical Physics, Trieste, Italy, where 
some years ago a stage of this investigation was dealt with. Inspiring discussions with 
Vo Thanh Cuong and Nguyen thi Hong Van are also acknowledged. This work was 
partially supported by the Vietnam National Research Program for Natural Sciences 
under Grant 411101.  
 
\end{document}